\documentclass[11pt]{amsart}

\usepackage{amssymb,mathrsfs,enumerate,hyperref,theoremref,graphicx,tikz,algorithm,algorithmic, tabulary,verbatim,array}
\usetikzlibrary{matrix,arrows}

\theoremstyle{mystyle}
\newtheorem{theorem}{Theorem}[section]

\DeclareMathOperator{\conv}{conv}

\DeclareMathOperator{\qdeg}{qdeg}
\DeclareMathOperator{\supp}{supp}
\DeclareMathOperator{\tdeg}{tdeg}

\newcommand\<{\langle}
\renewcommand\>{\rangle}

\newcommand\N{\mathbb{N}}
\newcommand\C{\mathbb{C}}
\newcommand\R{\mathbb{R}}
\newcommand\Q{\mathbb{Q}}
\newcommand\Z{\mathbb{Z}}
\renewcommand\k{\Bbbk}

\newcommand\bx{\mathbf{x}}

\newcommand\bu{\mathbf{u}}
\newcommand\bv{\mathbf{v}}

\newcommand\mm{\mathfrak{m}}
\newcommand\Ext{\text{Ext}}
\newcommand\Hom{\text{Hom}}
\newcommand\rank{{\rm rank}}

\newcommand\Oplus{\bigoplus}
\newcommand\varempty{\varnothing}
\newcommand\veps{\varepsilon}

\newcommand\vol{{\rm vol}}

\def\M2{\texttt{Macaulay2}}

\usepackage{xy}

\newcommand{\defn}[1]{\emph{#1}}

\begin{document}

\title{Computing quasidegrees of $A$-graded modules}
\author{Roberto Barrera}
\address[]{Department of Mathematics, Texas State University, San Marcos, Texas, 78666 U.S.A.}
\email{rbarrera@txstate.edu}
\date{}

\begin{abstract}
We describe the main functions of the \M2 package \texttt{Quasidegrees}.  The purpose of this package is to compute the quasidegree set of a finitely generated $\Z^d$-graded module presented as the cokernel of a monomial matrix.  We provide examples with motivation coming from $A$-hypergeometric systems.
\end{abstract}
\maketitle

\section{Introduction}\label{Intro}

Throughout $R=\k[x_1,\ldots,x_n]$ is a $\Z^d$-graded polynomial ring over a field $\k$ and $\mm=\< x_1,\ldots,x_n\>$ denotes the homogeneous maximal ideal in $R$.  Let $M=\Oplus_{\beta\in\Z^d} M_\beta$ be a $\Z^d$-graded $R$-module. The \defn{true degree set} of $M$ is
$$\tdeg(M)=\{\beta\in\Z^d\mid M_\beta\neq0\}.$$
The  \defn{quasidegree set} of $M$, denoted $\qdeg(M)$, is the Zariski closure in $\C^d$ of $\tdeg(M)$.  

The purpose of the \texttt{Macaulay2}\cite{M2} package \texttt{Quasidegrees}\cite{Qd} is to compute the quasidegree set of a finitely generated $\Z^d$-graded module presented as the cokernel of a monomial matrix.  By a monomial matrix, we mean a matrix where each entry is either zero or a monomial in $R$.  The initial motivation for \texttt{Quasidegrees} was to compute the quasidegree sets of certain local cohomology modules supported at $\mm$ of $\Z^d$-graded $R$-modules so there are some methods in the package specific to local cohomology.  Recall that the  \defn{ith local cohomology module} of $M$ with support at the ideal $I\subset R$ is the $i$th right derived functor of the left exact $I$-torsion functor
$$\Gamma_I(M)=\{m\in M\mid I^tm=0 \text{ for some }t\in\N\}$$
on the category of $R$-modules.

By the vanishing theorems of local cohomology \cite{CAAG}, the quasidegree sets of the local cohomology modules supported at $\mm$ of $M$ can be seen as measuring how far the module is from being Cohen-Macaulay.  From the $A$-hypergeometric systems point of view, the quasidegree set of the non-top local cohomology modules supported at $\mathfrak{m}$ of $R/I_A$, where $I_A$ is the toric ideal associated to $A$ in $R$, determine the parameters $\beta$ where the $A$-hypergeometric system $H_A(\beta)$ has rank higher than expected (see Section \ref{hypergeometricsystems}).

\section{Quasidegrees}\label{Quasidegrees}

The main function of \texttt{Quasidegrees} is \texttt{quasidegrees}, which computes the quasidegree set of a module that is presented by a monomial matrix.

We use the idea of standard pairs of monomial ideals to compute the quasidegree set of a $\Z^d$-graded $R$-module.  Given a monomial $x^\bu$ and a subset $Z\subset\{x_1,\ldots,x_n\}$, the pair $(x^\bu,Z)$ indexes the monomials $x^\bu\cdot x^\bv$ where $\supp(x^\bv)\subset Z$.  A \defn{standard pair} of a monomial ideal $I\subset R$ is a pair $(x^\bu,Z)$ satisfying:
\begin{enumerate}
\item $\supp(x^\bu)\cap Z=\varempty$,
\item all of the monomials indexed by $(x^\bu,Z)$ are outside of $I$,
\item $(x^\bu,Z)$ is maximal in the sense that $(x^\bu,Z)\nsubseteq(x^\bv,Y)$ for any other pair $(x^\bv,Y)$ satisfying the first two conditions.
\end{enumerate}
To compute the quasidegree set of $M$ we first find a monomial presentation of $M$ so that $M$ is the cokernel of a monomial matrix $\phi$.  We then compute the standard pairs of the ideals generated by the rows of $\phi$ and to each standard pair we associate the degrees of the corresponding variables.  The following algorithm is implemented in \texttt{Quasidegrees}.  The input is an $R$-module presented by a monomial matrix $\phi:R^s\rightarrow R^t$.  As in \M2, we write the degree of the $k$th factor of $R^t$ next to the $k$th row of the matrix $\phi$.

\begin{algorithm}[H]
\caption{Compute $\qdeg(M)$}
\label{alg1}
\begin{algorithmic}
\REQUIRE $R$-module $M$ presented by monomial matrix \\
	\hskip1cm$\phi=\alpha_i[c_{j,k}\bx^{\bu_{j,k}}]:R^s\rightarrow R^t$
\ENSURE $\qdeg(M)$
\STATE $Q=\varempty$
\FOR{ $1\leq k\leq t$} 
\STATE $SP=\{\text{standard pairs of  }\<c_{k,1}\bx^{\bu_{k,1}},c_{k,2}\bx^{\bu_{k,2}},\ldots,c_{k,s}\bx^{\bu_{k,s}}\>\}$
\STATE $Q=Q\cup\{\deg(\bx^\bu)+\alpha_k+\sum_{x_i\in F}\C\cdot\deg(x_i)\mid(\bx^\bu,Z)\in SP\}$
\ENDFOR
\RETURN Q
\end{algorithmic}
\end{algorithm}

In the implementation of Algorithm \ref{alg1} in \M2, we represent the output as a list of pairs $(\bu,Z)$ with $\bu\in\Q^d$ and $Z\subset\Q^d$ where the pair $(\bu,Z)$ represents the plane
$$\bu+\sum_{\bv\in Z}\C\cdot\bv.$$
The union of these planes over all such pairs in the output is the quasidegree set of $M$.

The following is an example of \texttt{Quasidegrees} computing the quasidegree set of an $R$-module:

\begin{verbatim}
i1 : R=QQ[x,y,z]
o1 = R
o1 : PolynomialRing
i2 : I=ideal(x*y,y*z)
o2 = ideal (x*y, y*z)
o2 : Ideal of R
i3 : M=R^1/I
o3 = cokernel | xy yz |
                             1
o3 : R-module, quotient of R
i4 : Q = quasidegrees M
o4 = {{0, {| 1 |}}, {0, {| 1 |, | 1 |}}}
o4 : List
\end{verbatim}

The above example displays a caveat of \texttt{quasidegrees} in that there may be some redundancies in the output. By a redundacy, we mean when one plane in the output is contained in another.  The redundancy above is clear:
$$\qdeg(\k[x,y,z]/\langle xy,yz\rangle)=\C=\{z_1+z_2\in\C\mid z_1,z_2\in\C\}.$$
The function \texttt{removeRedundancy} gets rid of redundancies in the list of planes:
\begin{verbatim}
i5 : removeRedundancy Q
o5 = {{0, {| 1 |, | 1 |}}}
o5 : List
\end{verbatim}

\section{Quasidegrees and hypergeometric systems}\label{hypergeometricsystems}%

In this section, we discuss the motivation for \texttt{Quasidegrees} and the methods in \texttt{Quasidegrees} that aid us in our studies. Let $A=[a_1~a_2~\cdots~a_n]$ be an integer $(d\times n)$-matrix with $\Z A=\Z^d$ and such that the cone over its columns is pointed.  There is a natural $\Z^d$-grading of $R$ by the columns of $A$ given by $\deg(x_j)=a_j$, the $j$th column of $A$.  A module that is homogeneous with respect to this grading is said to be {\it A-graded}.  By the assumptions on $A$, $R$ is positively graded by $A$, that is, the only polynomials of degree 0 are the constants.  Given such a matrix $A$ and a polynomial ring $R$ in $n$ variables, the method \texttt{toGradedRing} gives $R$ an $A$-grading.  For example, let $A=\left(\begin{smallmatrix}1&1&1&1&~1\\0&0&1&1&~0\\0&1&1&0&-2\end{smallmatrix}\right)$.  We make the $A$-graded polynomial ring $\Q[x_1,x_2,x_3,x_4,x_5]$ :
 
\begin{verbatim}
i6 : A=matrix{{1,1,1,1,1},{0,0,1,1,0},{0,1,1,0,-2}}
o6 = | 1 1 1 1 1  |
     | 0 0 1 1 0  |
     | 0 1 1 0 -2 |
              3        5
o6 : Matrix ZZ  <--- ZZ
i7 : R=QQ[x_1..x_5]
o7 = R
o7 : PolynomialRing
i8 : R=toGradedRing(A,R)
o8 = R
o8 : PolynomialRing
i9 : describe R
o9 = QQ[x , x , x , x , x , Degrees => {{1}, {1}, {1}, {1},
         1   2   3   4   5              {0}  {0}  {1}  {1} 
                                        {0}  {1}  {1}  {0} 
     -------------------------------------------------------
     {1 }}, Heft => {1, 2:0}, MonomialOrder =>
     {0 }                                     
     {-2}                                     
     -------------------------------------------------------
     {MonomialSize => 32}, DegreeRank => 3]
     {GRevLex => {5:1}  }
     {Position => Up    }

\end{verbatim}

The \defn{toric ideal associated to} $A$ in $R$ is the binomial ideal
$$I_A=\langle\bx^\bu-\bx^\bv:A\bu=A\bv\rangle.$$
The method \texttt{toricIdeal} computes the toric ideal associated to $A$ in the ring $R$.  We continue with the $A$ and $R$ from the above example and compute the toric ideal $I_A$ associated to $A$ in $R$:

\begin{verbatim}
i10 : I=toricIdeal(A,R)
                            2    2     2              3  
o10 = ideal (x x  - x x , x x  - x x , x x  - x x x , x  -
             1 3    2 4   1 4    3 5   1 4    2 3 5   1  
     -------------------------------------------------------
      2
     x x )
      2 5

o10 : Ideal of R

\end{verbatim}

We now introduce $A$-hypergeometric systems.  Given a matrix $A\in\Z^{d\times n}$ as above and a $\beta\in\C^d$, the \defn{A-hypergeometric system with parameter $\beta\in\C^d$} \cite{SST}, denoted $H_A(\beta)$,  is the system of partial differential equations:
\begin{align*}
&\frac{\partial^{|\bv|}}{\partial\bx^\bv}\phi(\bx)=\frac{\partial^{|\bu|}}{\partial\bx^\bu}\phi(\bx)\text{ for all $\bu,\bv$, $A\bu=A\bv$}\\
&\sum_{j=1}^n a_{ij}x_j\frac{\partial}{\partial x_j}\phi(\bx)=\beta_i\phi(\bx),\text{ for }i=1,\ldots, d.
\end{align*}

Such systems are sometimes called \defn{GKZ-hypergeometric systems}.  The function \texttt{gkz} in the \texttt{Macaulay2} package \texttt{Dmodules} computes this system as an ideal in the Weyl algebra.  The \defn{rank} of $H_A(\beta)$ is
{\small
$$\rank(H_A(\beta))=\dim_\C\left\{\begin{array}{l}\text{germs of holomorphic solutions of $H_A(\beta)$}\\\text{near a generic nonsingular point}\end{array}\right\}.$$}

The function \texttt{holonomicRank} in \texttt{Dmodules} computes the rank of an $A$-hypergeometric system.  In general, rank is not a constant function of $\beta$. Denote $\vol(A)$ to be $d!$ times the Euclidean volume of $\conv(A\cup\{0\})$ the convex hull of the columns of $A$ and the origin in $\R^d$.  The following theorem gives the parameters $\beta$ for which $\rank(H_A(\beta))$ is higher than expected:

\begin{theorem}\label{rankthm}\cite{MMW}
Let $H_A(\beta)$ be an $A$-hypergeometric system with parameter $\beta$.  If $\beta\in\qdeg(\bigoplus_{i=0}^{d-1}H_\mm^i(R/I_A))$ then $\rank(H_A(\beta))>\vol(A)$.  Otherwise, $\rank(H_A(\beta))=\vol(A)$.
\end{theorem}

Since Theorem \ref{rankthm} was the initial motivation for \texttt{Quasidegrees}, the package has a method \texttt{quasidegreesLocalCohomology} (abbreviated \texttt{qlc}) to compute the quasidegree set of the local cohomology modules $H_\mm^i(R/I_A)$.  If the input is an integer $i$ and the $R$-module $R/I_A$, then the method computes $\qdeg(H^i_\mm(R/I_A))$.  If the input is only the module $R/I_A$, the method computes the quasidegree set in Theorem \ref{rankthm}.

We use graded local duality to compute the local cohomology modules of a finitely generated $A$-graded $R$-module supported at the maximal ideal $\mm$:

\begin{theorem}\label{GLC}(Graded local duality \cite{BH98,M01})
Given an $A$-graded $R$-module $M$, there is an $A$-graded vector space isomorphism
$$\Ext_R^{n-i}(M,R)_\alpha\cong \Hom_\k(H_\mm^i(M)_{-\alpha-\veps_A},\k)$$
where $\mm=\langle x_1,\ldots,x_n\rangle$ and $\veps_A=\sum_{j=1}^n a_j$.
\end{theorem}

The algorithm implemented for \texttt{quasidegreesLocalCohomology} is essentially Algorithm \ref{alg1} applied to the $\Ext$-modules of $M$ with the additional twist of $\varepsilon_A$ coming from local duality.  For our purposes, we exploit the fact that the higher syzygies of $R/I_A$ are generated by monomials in $R^m$ (see \cite{CCA}, Chapter 9).

Continuing our running example, we use \texttt{quasidegreesLocalCohomology} to compute the quasidegree set of $\bigoplus_{i=0}^{d-1}H_\mm^i(R/I_A)$:

\begin{verbatim}
i11 : M=R^1/I
o11 = cokernel | x_1x_3-x_2x_4 x_1x_4^2-x_3^2x_5 
x_1^2x_4-x_2x_3x_5 x_1^3-x_2^2x_5 |
                            1
o11 : R-module, quotient of R
i12 : quasidegreesLocalCohomology M
o12 = {{| 0 |, {| 1  |}}}
        | 0 |   | 0  |
        | 1 |   | -2 |
o12 : List
\end{verbatim}

Thus
\begin{equation}\label{plane}
\qdeg\left(\bigoplus_{i=0}^{d-1}H_\mm^i(R/I_A)\right)=\left[\begin{smallmatrix}0\\0\\1\end{smallmatrix}\right]+\C\cdot\left[\begin{smallmatrix}\phantom{-}1\\\phantom{-}0\\-2\end{smallmatrix}\right].\end{equation}

As a check, we use the methods \texttt{gkz} and \texttt{holonomicRank} from the package \texttt{Dmodules} to compute $\rank(H_A(0))$ and $\rank(H_A(\beta))$ for two different $\beta$ in (\ref{plane}) and demonstrate a rank jump:

\begin{verbatim}
i13 : holonomicRank gkz(A,{0,0,0}) -- vol A in this case
o13 = 4
i14 : holonomicRank gkz(A,{0,0,1})
o14 = 5
i15 : holonomicRank gkz(A,{3/2,0,-2})
o15 = 5
\end{verbatim}

\bibliographystyle{alpha}
\bibliography{bib1}

\end{document}